\documentclass[12pt]{amsart}

\usepackage{amsmath,amssymb,amsbsy,amsfonts,amsthm,latexsym,
                        amsopn,amstext,amsxtra,euscript,amscd,mathrsfs,color,bm,cite}

\usepackage{float}
\usepackage[english]{babel}
\usepackage{mathtools}
\usepackage{todonotes}
\usepackage{url}
\usepackage[colorlinks,linkcolor=blue,anchorcolor=blue,citecolor=blue,backref=page]{hyperref}

\usepackage[norefs,nocites]{refcheck}

\newtheorem{theorem}{Theorem}
\newtheorem{lemma}[theorem]{Lemma}

\numberwithin{equation}{section}
\numberwithin{theorem}{section}
\numberwithin{table}{section}

\numberwithin{figure}{section}

\newfont{\teneufm}{eufm10}
\newfont{\seveneufm}{eufm7}
\newfont{\fiveeufm}{eufm5}
\newfam\eufmfam
                \textfont\eufmfam=\teneufm \scriptfont\eufmfam=\seveneufm
                \scriptscriptfont\eufmfam=\fiveeufm
\def\bbbc{{\mathchoice {\setbox0=\hbox{$\displaystyle\rm C$}\hbox{\hbox
to0pt{\kern0.4\wd0\vrule height0.9\ht0\hss}\box0}}
{\setbox0=\hbox{$\textstyle\rm C$}\hbox{\hbox
to0pt{\kern0.4\wd0\vrule height0.9\ht0\hss}\box0}}
{\setbox0=\hbox{$\scriptstyle\rm C$}\hbox{\hbox
to0pt{\kern0.4\wd0\vrule height0.9\ht0\hss}\box0}}
{\setbox0=\hbox{$\scriptscriptstyle\rm C$}\hbox{\hbox
to0pt{\kern0.4\wd0\vrule height0.9\ht0\hss}\box0}}}}
\def\bbbq{{\mathchoice {\setbox0=\hbox{$\displaystyle\rm
Q$}\hbox{\raise 0.15\ht0\hbox to0pt{\kern0.4\wd0\vrule
height0.8\ht0\hss}\box0}} {\setbox0=\hbox{$\textstyle\rm
Q$}\hbox{\raise 0.15\ht0\hbox to0pt{\kern0.4\wd0\vrule
height0.8\ht0\hss}\box0}} {\setbox0=\hbox{$\scriptstyle\rm
Q$}\hbox{\raise 0.15\ht0\hbox to0pt{\kern0.4\wd0\vrule
height0.7\ht0\hss}\box0}} {\setbox0=\hbox{$\scriptscriptstyle\rm
Q$}\hbox{\raise 0.15\ht0\hbox to0pt{\kern0.4\wd0\vrule
height0.7\ht0\hss}\box0}}}}
\def\bbbt{{\mathchoice {\setbox0=\hbox{$\displaystyle\rm
T$}\hbox{\hbox to0pt{\kern0.3\wd0\vrule height0.9\ht0\hss}\box0}}
{\setbox0=\hbox{$\textstyle\rm T$}\hbox{\hbox
to0pt{\kern0.3\wd0\vrule height0.9\ht0\hss}\box0}}
{\setbox0=\hbox{$\scriptstyle\rm T$}\hbox{\hbox
to0pt{\kern0.3\wd0\vrule height0.9\ht0\hss}\box0}}
{\setbox0=\hbox{$\scriptscriptstyle\rm T$}\hbox{\hbox
to0pt{\kern0.3\wd0\vrule height0.9\ht0\hss}\box0}}}}
\def\bbbs{{\mathchoice
{\setbox0=\hbox{$\displaystyle     \rm S$}\hbox{\raise0.5\ht0\hbox
to0pt{\kern0.35\wd0\vrule height0.45\ht0\hss}\hbox
to0pt{\kern0.55\wd0\vrule height0.5\ht0\hss}\box0}}
{\setbox0=\hbox{$\textstyle        \rm S$}\hbox{\raise0.5\ht0\hbox
to0pt{\kern0.35\wd0\vrule height0.45\ht0\hss}\hbox
to0pt{\kern0.55\wd0\vrule height0.5\ht0\hss}\box0}}
{\setbox0=\hbox{$\scriptstyle      \rm S$}\hbox{\raise0.5\ht0\hbox
to0pt{\kern0.35\wd0\vrule height0.45\ht0\hss}\raise0.05\ht0\hbox
to0pt{\kern0.5\wd0\vrule height0.45\ht0\hss}\box0}}
{\setbox0=\hbox{$\scriptscriptstyle\rm S$}\hbox{\raise0.5\ht0\hbox
to0pt{\kern0.4\wd0\vrule height0.45\ht0\hss}\raise0.05\ht0\hbox
to0pt{\kern0.55\wd0\vrule height0.45\ht0\hss}\box0}}}}
\def\bbbz{{\mathchoice {\hbox{$\sf\textstyle Z\kern-0.4em Z$}}
{\hbox{$\sf\textstyle Z\kern-0.4em Z$}} {\hbox{$\sf\scriptstyle
Z\kern-0.3em Z$}} {\hbox{$\sf\scriptscriptstyle Z\kern-0.2em
Z$}}}}

\def\squareforqed{\hbox{\rlap{$\sqcap$}$\sqcup$}}
\def\qed{\ifmmode\squareforqed\else{\unskip\nobreak\hfil
\penalty50\hskip1em\null\nobreak\hfil\squareforqed
\parfillskip=0pt\finalhyphendemerits=0\endgraf}\fi}

\def\cG{{\mathcal G}}

\def\cS{{\mathcal S}}

\def\le{\leqslant}
\def\leq{\leqslant}
\def\ge{\geqslant}
\def\leq{\leqslant}

\newcommand{\ignore}[1]{}

\def\e{\mathbf{e}}

\def\xbar{\overline{x}}
\def\sbar{\overline{s}}

\def \F{\mathbb{F}}

\def \Z{\mathbb{Z}}

\def \Z{\mathbb{Z}}

\def\mand{\qquad\mbox{and}\qquad}

\def\\{\cr}
\def\({\left(}
\def\){\right)}

\def\e{{\mathbf{\,e}}}
\def\ep{{\mathbf{\,e}}_p}

\newcommand{\lcm}{\operatorname{lcm}}

\begin{document}

\title[Modular inverses from short intervals]
{On the distribution of modular inverses from short intervals}

\author{Moubariz~Z.~Garaev}
\address{Centro  de Ciencias Matem{\'a}ticas,  Universidad Nacional Aut\'onoma de
M{\'e}\-xico, C.P. 58089, Morelia, Michoac{\'a}n, M{\'e}xico}
\email{garaev@matmor.unam.mx}

\author[I. E. Shparlinski] {Igor E. Shparlinski}
\address{School of Mathematics and Statistics, University of New South Wales, Sydney, NSW 2052, Australia}
\email{igor.shparlinski@unsw.edu.au}

\begin{abstract}  
For a prime number $p$ and integer $x$ with $\gcd(x,p)=1$ let $\xbar$ denote the multiplicative inverse
of $x$ modulo $p.$ In the present paper we are interested in the problem of distribution
modulo $p$ of the sequence
$$
\xbar, \qquad x =1, \ldots, N, 
$$
and in lower bound estimates for the corresponding exponential sums. As representative examples, we state the following two consequences of the main results. 

For  any fixed  $A > 1$ and for any sufficiently large integer $N$ there exists a prime number $p$ with
$$
(\log p)^A \asymp N
$$
such that   
$$
\max_{(a,p)=1}\left|\sum_{x\le N}\e_p(a\xbar)\right|\gg N. 
$$

For any fixed positive  $\gamma< 1$  
there exists a positive constant $c$  such that the following holds: for any sufficiently large integer $N$ there is a prime number $p > N$ such that  
$$
N > \exp\(c(\log p\log\log p)^{\gamma/(1+\gamma)}\)
$$
and 
$$
\max_{(a,p)=1}\left|\sum_{x\le N}\e_p\(a \xbar\)\right|\gg N^{1-\gamma}.
$$
\end{abstract}

\keywords{modular inverses, uniform distribution. dicrepancy}
\subjclass[2020]{11J71, 11K38,  11L05,  11N25}	

\maketitle

\tableofcontents
\section{Introduction}

\subsection{Motivation} 
Given a prime $p$ and a positive integer $N < p$, we consider 
the sequence of modular inverses 
\begin{equation}
\label{eq: x-bar}
\xbar, \qquad x =1, \ldots, N,
\end{equation}
defined by the conditions
$$
\xbar x \equiv 1 \pmod p \mand 1 \le \xbar < p. 
$$

There is an extensive literature dedicated to studying the distribution 
of the sequence~\eqref{eq: x-bar}, with application to an amazingly large variety
of problems in number theory and beyond, see the survey~\cite{Shp} and also more recent works~\cite{Baier,BrHay,Chan,Hump,Ust}.

When $N \ge p^{1/2+\varepsilon}$, the classical Weil bound of Kloosterman sums,
see~\cite[Theorem~11.11]{IwKow}, immediately implies that the elements of the sequence~\eqref{eq: x-bar},
are sufficiently uniformly distributed in interval $[1, p-1]$. For smaller values of $N< p^{1/2}$ the problem becomes 
much more difficult and only significantly  weaker results are known, see, for example,~\cite{BG,Kar1,Kar2, Kor, Kor1, Kor2}. 
It is known from the work of Bourgain and Garaev~\cite{BG} that one has the  uniformity of distribution in the range
\begin{equation}
\label{eq: NontrivBound}
p^{1/2}\ge N \ge \exp\(\(\log p\)^{2/3}(\log\log p)^{c} \), 
\end{equation}
where  $c>3$ is any fixed constant. The constant $2/3$ that appears in  $(\log p)^{2/3}$ is the best known up to the date, while the constant $c$ is claimed to be improved in Korolev~\cite[Theorem~2]{Kor2} (the reader should be aware, however, that the proof of~\cite[Lemma~2]{Kor2} contains a gap, 
while the results of Snurnitsyn~\cite{Snur} are, unfortunately, incorrect). In the range of $N$ just below $p^{1/2}$, 
and a nontrivial bound is given by Fouvry,   Kowalski,
 Michel,  Raju,  Rivat an  Soundararajan~\cite[Corollary~1.7]{FKMRRS}.

Here we approach the question from a different direction and ask what are the largest possible values of $N$ 
for which the sequence~\eqref{eq: x-bar} is not uniformly distributed or is rather poorly distributed. 

\subsection{Formal set-up and notation} 
As usual, we measure the quality of the distribution of the sequence~\eqref{eq: x-bar}
via its discrepancy 
\begin{equation}
\label{eq: DpN}
D_p(N) = \sup_{\alpha \in [0,1]} \left| \frac{1}{N} A_p(N, \alpha) - \alpha  \right|
\end{equation}
where
$$A_p(N, \alpha) =\#\left\{x\in \{1, \ldots, N\}:~ \xbar/p \le \alpha\right\}.
$$
and  for a finite set $\cS$ we use $\# \cS$ to denote its cardinality. 

We say that  the sequence~\eqref{eq: x-bar} is uniformly distributed if $D_p(N) \to 0$ as $p,N\to \infty$.

We recall
that the notations $U = O(V)$,  $U \ll V$ and  $V \gg U$  are
all equivalent to the statement that $|U| \le c V$ holds
with some constant $c> 0$.

Any implied constants in symbols $O$, $\ll$
and $\gg$ may occasionally, where obvious, depend on the  real positive parameters $A$ or $\gamma$, 
 and are absolute otherwise. 
 We also write 
 $$U \asymp V
 $$ 
 as an equivalent of  $U\ll V \ll U$.   

Finally $U= o(V)$ means that $U \le \psi(V)V$ for some
function $\psi$ such that $\psi(V) \to 0$ as $V \to \infty$.

We also write $\e(z)=\exp(2\pi i z)$ and $\ep(z)=\e(z/p)$.

\subsection{Main results}

First we consider the case of non-uniform distribution on rather long segments~\eqref{eq: x-bar}. 
\begin{theorem}
\label{thm: Non u.d.}  For  any fixed  $A > 1$ and any sufficiently large positive integer $N$ there exists a  prime number $p$ with
$$
(\log p)^A \asymp N 
$$
such that  $D_p(N) \gg 1$. 
\end{theorem}

Using~\cite[Chapter~1,
Theorem~1]{Mont} which is  a version of the celebrated Erd\H{o}s--Tur{\'a}n--Koksma inequality (see, for example,~\cite[Theorem~1.21]{DrTi}),   given  Lemma~\ref{lem:Mont} below), which bound the discrepancy in terms 
of exponential sums, we obtain the following  estimate on certain $L_1$-moment of exponential sums 
with reciprocals.

\begin{theorem}
\label{thm: Exp Sum non saving}  For  any fixed  $A > 1$ there exists an integer $k_0>0$ such that the following holds: for any sufficiently large integer $N$ there exists a prime number $p$ with
$$
(\log p)^A \asymp N
$$
such that   
$$
\sum_{a=1}^k\left|\sum_{x\le N}\e_p(a\xbar)\right|\gg kN, 
$$ 
for any $k_0\le k\le N.$
\end{theorem}

Next we consider the case when the discrepancy may be nontrivial but admits only
a  small power saving. 

\begin{theorem}
\label{thm: Discrep power non saving} 
For any fixed positive  $\gamma< 1$ there exist positive constants $c_1$ and $c_2$ such that
the following holds: for any sufficiently large integer N there exists a prime number $p$
such that   
$$
\exp\(c_1(\log p\log\log p)^{\gamma/(1+\gamma)}\) > N > \exp\(c_2(\log p\log\log p)^{\gamma/(1+\gamma)}\)
$$
and 
$$
D_p(N)\gg  N^{-\gamma} .
$$
\end{theorem}

We  then derive an analogue of  Theorem~\ref{thm: Exp Sum non saving}.

\begin{theorem}
\label{thm: Exp Sum non large power saving} 
For any fixed positive  $\gamma< 1$  there exist positive constants $c_1$, $c_2$ and $\beta$ 
such that the following holds: for any sufficiently large integer $N$ there is a 
prime number $p$ such that  
$$
\exp\(c_1(\log p\log\log p)^{\gamma/(1+\gamma)}\) > N > \exp\(c_2(\log p\log\log p)^{\gamma/(1+\gamma)}\)
$$
and 
$$
\sum_{a=1}^k\left|\sum_{x\le N}\e_p\(a \xbar\)\right|\gg kN^{1-\gamma},
$$
for any $\beta N^{\gamma}<k<N$.
\end{theorem}

In particular, we immediately see that under the conditions of Theorems~\ref{thm: Exp Sum non saving}  
and~\ref{thm: Exp Sum non large power saving}, we have 
\begin{equation}
\label{eq: Max Exp Sum}
\max_{\gcd(a,p)=1} \left|\sum_{x\le N}\e_p\(a \xbar\)\right|\gg N 
\end{equation}
and 
$$
\max_{\gcd(a,p)=1} \left|\sum_{x\le N}\e_p\(a \xbar\)\right|\gg N^{1-\gamma},
$$
 respectively.  The case of $\gamma$ approaching $1/2$ from below  is especially 
 interesting. In particular,  Theorem~\ref{thm: Exp Sum non large power saving} means that 
 non square-root cancellation in exponential sums is possible for 
 $$
N \le  \exp\((\log p)^{1/3}\).
$$
 
\section{Preliminaries}

\subsection{Prime numbers}

First we recall the famous result of Linnik about primes in an arithmetic progression.
We present it in the currently strongest form due to Xylouris~\cite{Xy}.

\begin{lemma}
\label{lem:Linnik} For any integers $m > a \ge 1$ with $\gcd(a,m)=1$, there is a prime 
$p \equiv a \pmod m$ with $p \ll m^{5.18}$. 
\end{lemma} 

We remark that in our applications, the modulus $q$ is a product of small primes, so
we can replace the power $5.18$ in Lemma~\ref{lem:Linnik} with  $12/5+o(1)$ via a result of Chang~\cite{Chang}. 

Below we  use Lemma~\ref{lem:Linnik} with $a =m-1$, thus we have $m-1 \le p  \ll m^{5.18}$
which we use in a crude form 
$$
\log p \asymp \log m.
$$

\subsection{Smooth numbers}
Let $\Psi(x, y)$ denote the number of $y$-smooth numbers less than or equal to $x$, that is, 
$\Psi(x, y)$ denote the number of  integers up to  $x$ without prime divisors exceeding $y$, we refer to~\cite{Gran,HilTen} 
exhaustive outline of most important facts on smooth numbers.

Let $\rho(u)$  be the {\it Dickman--de Bruijn function} which  is defined recursively by:
$$\rho(u)= 
\begin{cases} 
1 & \text{ if } 0 \leq u \leq 1,\\
{\displaystyle 1 - \int_{n}^{t}\frac{\rho(v-1)}{v}dv} & \text{ if } u>1.
\end{cases}
$$

We make use of the following explicit lower bound on $\Psi(x, y)$,  which is
due to S.~Konyagin and C. Pomerance~\cite[Theorem~2.1]{KoPo}.

\begin{lemma}
\label{lem:KPsmooth} If $x\ge 4$ and $x \ge y \ge 2$, then
$$
\Psi(x,y) > x^{1-\log \log x/\log y}.
$$
\end{lemma}

We recall the following asymptotic formula, see~\cite[Part~III, Corollary~5.19]{Ten}

\begin{lemma}
\label{lem:Smooth} Let $\varepsilon > 0$ be fixed.  If 
$$
x\ge 4 \mand x \ge y \ge \exp\(\(\log \log x\)^{5/3+\varepsilon}\),
$$ 
then
$$
\Psi(x,y) = \rho(u) x\(1+ O\(\frac{\log(u+1)}{\log y}\)\), 
$$
where 
$$
u= \frac{\log x}{\log y}.
$$ 
\end{lemma}

Finally,  we have the following result obtained, after simple manipulations, by combining~\cite[Theorem~3]{HilTen1} 
and~\cite[Equation~(2.4)]{HilTen1}:

\begin{lemma}
\label{lem:psicxHildebrand}
For $x\geq y\geq 2$ and $1\leq c\leq y$, we have  
\begin{align*}
    \Psi(cx, y)=\Psi(x, y)& \left(1+\frac{y}{\log x}\right)^{\log c/\log y}\\
    &  \quad \times \left(1+O\left(\frac{\log\log(1+y)}{\log^2 y}\cdot\log\left(1+\frac{y}{\log x}\right)\right)\right)\\
    & \qquad  \qquad \qquad    \qquad \qquad  \ \times \left(1+O\left(\frac{1}{u}+\frac{\log y}{y}\right)\right), 
\end{align*}
where
$$
u= \frac{\log x}{\log y}.
$$ 
\end{lemma}

 \subsection{Discrepancy and exponential sums}
 We also need the
following result, which is a particular case of~\cite[Chapter~1,
Theorem~1]{Mont}.

\begin{lemma}
\label{lem:Mont} Let $\gamma_1, \ldots, \gamma_N$ be a sequence
of $N$ points of the interval $[0,1]$. Then for any integer
$k\ge 1$ and any  $\alpha\in [0,1]$, we
have
\begin{align*}
\left|\# \{n =1, \ldots, N:~\gamma_n  \in [0, \alpha]\} - \alpha N 
\right| &\\
\le \frac{N}{k+1} +  \(\frac{2}{k+1}+2\alpha\) \sum_{a=1}^k &
\left|\sum_{n=1}^N \exp(2 \pi i a \gamma_n)\right|.
\end{align*}
\end{lemma}

\section{Proof of  results about trivial bounds}

\subsection{Proof of  Theorem~\ref{thm: Non u.d.}} 
Let us fix some real  $A\ge 1$. It is more convenient to get a lower bound on  $D_p(2N)$, 
which,  of course is an equivalent question.

For a sufficiently large integer $N$, let $\cS$ be the set  of $(2N)^{1/A}$-smooth numbers in the interval $[N, 2N]$. 
Then, by Lemma~\ref{lem:Smooth}, we have 
\begin{equation}
\label{eq: Set S}
\# \cS =  \rho(A) N + o(N)  \ge 0.5  \rho(A) N , 
\end{equation}
provided that $N$ is large enough. 

Next, let 
$$
m = \lcm[s:~s \in \cS]
$$
be the least common multiple of the elements of $\cS$.  By the prime number theorem $m$ is a 
product of at most 
$$
\pi\((2N)^{1/A}\) \ll N^{1/A}/\log N
$$ 
prime powers $q^r$ with a prime $q\le (2N)^{1/A}$ and such that $q^r\le 2N$.  
Therefore, 
$$
m \le (2N)^{O\(N^{1/A}/\log N\)} =  \exp\(O\(N^{1/A}\)\).
$$
One can easily see there is a matching lower bound on $m$, thus 
$$
\log m \asymp  N^{1/A}.
$$
Therefore, by  Lemma~\ref{lem:Linnik} there is a prime number $p$ with 
$$
p \equiv - 1  \pmod {m}  \mand  \log p  \asymp \log m  \asymp N^{1/A}. 
$$

We now take 
\begin{equation}
\label{eq: alpha0}
\alpha_0 = \frac{2}{N}
\end{equation}
in the definition of the discrepancy $D_p(2N)$ in~\eqref{eq: DpN}. 

Next we observe that for the above choice of $p$, for every $s\in S$ we have  
$\sbar  = (p+1)/s$ and thus 
\begin{equation}
\label{eq: small sbar}
\sbar/p \le (1+1/p)/s \le 2/s \le 2/N = \alpha_0.
\end{equation}
This, together with~\eqref{eq: Set S} implies that 
\begin{equation}
\label{eq: A(alpha0)}
A_p(2N, \alpha_0)  \ge \# \cS \ge 0.5  \rho(A) N
\end{equation}
and hence 
$$
 \frac{1}{2N}  A_p(2N, \alpha_0) - \alpha_0 \ge  0.25  \rho(A)  - 2/N \ge 0.2   \rho(A)
$$
provided that $N$ is large enough and hence  $D_p(2N) \gg 1$.

\subsection{Proof of Theorem~\ref{thm: Exp Sum non saving}}
We continue to work with $2N$ instead of $N$ and we use the same construction of $p$ and $N$ as in Theorem~\ref{thm: Non u.d.}

We also choose $\alpha_0$ as in~\eqref{eq: alpha0} so we have~\eqref{eq: A(alpha0)}.

We now assume that $k\ge  8/\rho(A)$, so that  we see from~\eqref{eq: A(alpha0)} that 
$$
2N/(k+1) \le 0.5  A_p(2N, \alpha_0).
$$ 
Then from Lemma~\ref{lem:Mont} and  using~\eqref{eq: A(alpha0)} again, we conclude that 
$$
  \(\frac{2}{k+1}+2\alpha_0\)  \sum_{a=1}^k\left|\sum_{x\le 2N}\e_p(a\xbar)\right|\ge 
   0.5  A_p(2N, \alpha_0)  \ge 0.25  \rho(A) N.
$$
From the choice of $\alpha_0$  in~\eqref{eq: alpha0}, for  $k \le N$, we conclude that
$$
  \frac{2}{k+1}+2\alpha_0 \le   \frac{6}{k+1}
$$ 
and the result follows.

\section{Proof of  results about bounds with small polynomial saving}

\subsection{Proof of  Theorem~\ref{thm: Discrep power non saving}}
For a sufficiently large integer $N$, let $\cS$ be the set  of $(\log N)^{1/\gamma}$-smooth numbers in the interval $[N, 2N]$. 
Then, by Lemmas~\ref{lem:KPsmooth} and~\ref{lem:psicxHildebrand}, we have 
\begin{equation}
\label{eq: again Set S}
\# \cS\gg  N^{1-\gamma} , 
\end{equation}
provided that $N$ is large enough. 

As in the proof of Theorem~\ref{thm: Non u.d.}, we let  
$$
m = \lcm[s:~s \in \cS]
$$
be the least common multiple of the elements of $\cS$.  By the prime number theorem $m$ is a 
product of at most 
$$
\pi\((\log N)^{1/\gamma}\) \ll (\log N)^{1/\gamma}/\log\log N
$$
prime powers $q^r$ with a prime $q\le (\log N)^{1/\gamma}$ and such that $q^r\le 2N$. 
Therefore, 
$$
m \le (2N)^{O\((\log N)^{1/\gamma}/(\log \log N)\)} =  \exp\(O\(\(\log N\)^{1+1/\gamma}/\log\log N\)\)
$$
and as in the proof of Theorem~\ref{thm: Non u.d.} we note that there is matching lower bound on $m$, thus 
$$
 m   \asymp \(\log N\)^{1+1/\gamma}/\log\log N
$$  
Therefore, by  Lemma~\ref{lem:Linnik} there is a prime number $p$ with 
$$
p \equiv - 1 \pmod m \mand  \log p \asymp m  \asymp \(\log N\)^{1+1/\gamma}/\log\log N. 
$$
that is,
$$
\exp\(c_1(\log p\log\log p)^{\gamma/(1+\gamma)}\) > N > \exp\(c_2(\log p\log\log p)^{\gamma/(1+\gamma)}\),
$$
for some constants $c_1\ge c_2 >0$. 
We now take 
$$
\alpha_0 = \frac{2}{N}
$$
in the definition of the discrepancy $D_p(2N)$ in~\eqref{eq: DpN} and observe 
that we still have~\eqref{eq: small sbar}.

This, together with~\eqref{eq: again Set S} applies that 
$$
A_p(2N, \alpha_0)  \ge \# \cS \gg N^{1-1/\gamma}
$$
and hence 
\begin{equation}
\label{eqn: A(alpha0)/2N -alpha0}
 \frac{1}{2N}  A_p(2N, \alpha_0) - \alpha_0 \gg N^{-1/\gamma}.
\end{equation}
provided that $N$ is large enough.

\subsection{Proof of Theorem~\ref{thm: Exp Sum non large power saving}}
We use the set $\cS$ from the proof of Theorem~\ref{thm: Discrep power non saving}
and then proceed as in the proof of Theorem~\ref{thm: Exp Sum non saving},
using~\eqref{eqn: A(alpha0)/2N -alpha0}.

\section{Comments} 

We note that the non-uniformity of distribution result of Theorem~\ref{thm: Non u.d.}  
as well  the lower bound~\eqref{eq: Max Exp Sum} on exponential sums, are in stark contrast 
with what is expected for the residues of small multiplicative subgroups of $\F_p^*$.
Namely,   the conjecture of Montgomery,   Vaughan and  Wooley~\cite{MontVW}
suggests that the residues of the multiplicative subgroups   $\cG \subseteq \F_p^*$ are 
uniformly distributed modulo $p$, provided $\#\cG/\log p \to \infty$ as $p\to \infty$, see also~\cite{Bour1, Bour2}
for some  results  towards this conjecture. 

On the other hand, nontrivial upper bounds are also known for shorter segments 
of reciprocals, see~\eqref{eq: NontrivBound}, rather than for subgroups   $\cG \subseteq \F_p^*$, 
where one currently requires 
$$
\# \cG \ge \exp\(\frac{\log p}{\log \log p}\)
$$
with some constant $c>0$, see~\cite{Bour2}.

It is easy to see that   Theorems~\ref{thm: Non u.d.}  and~\ref{thm: Exp Sum non saving} 
can be extended to the discrepancy  and exponential sums associated with Laurent polynomials
$$
 \sum_{i=-u}^v a_i X^i \in \Z[X]
$$
with fixed integer coefficients.  It has been known since the works of Bourgain~\cite{Bour3} and
Fouvry~\cite{Fou} that short exponential sums with $a\xbar^2$ are important for studying 
fundamental solution to some families of Pell equations. Such sums also appear in 
studying the distribution of square-free numbers in arithmetic progressions, see~\cite{Nun}.

Finally, we note that our approach also allows us to get results, which ``interpolate" between 
the  trivial bounds of  Theorems~\ref{thm: Non u.d.}  and~\ref{thm: Exp Sum non saving} 
 the bounds with a power saving 
of Theorems~\ref{thm: Discrep power non saving} and~\ref{thm: Exp Sum non large power saving}. 
For example, one can consider the case 
when $D_p(N)$ admits only a logarithmic saving, that is, $D_p(N) \gg (\log N)^{-\gamma}$ 
for some $\gamma >0$.

\section*{Acknowledgement}

The authors would like to thanks Zeev Rudnick, discussions with whom stimulated this work. 

During the preparation of this work, I.S. was supported in part by the Australian Research Council Grants~DP230100534.

\end{document}